\documentclass[12pt]{amsart}
\usepackage{amssymb,amsmath}
\makeatletter
\def\@cite#1#2{{\m@th\upshape\bfseries%
[{#1\if@tempswa{\m@th\upshape\mdseries, #2}\fi}]}}
\makeatother
\theoremstyle{plain}
\newtheorem{thm}{Theorem}[section]
\newtheorem{lem}[thm]{Lemma}
\newtheorem{cor}[thm]{Corollary}
\newtheorem{prop}[thm]{Proposition}
\theoremstyle{definition}
\newtheorem{rem}[thm]{Remark}

\newcommand{\Prf}{\noindent\textbf{Proof.\ }}
\newcommand{\bx}{\strut\hfill$\blacksquare$\medbreak}

\newcommand{\ca}{\mathrm{C}^*}

\newcommand{\ol}{\overline}

\newcommand{\bbA}{{\mathbb{A}}}

\newcommand{\bbF}{{\mathbb{F}}}

\newcommand{\bbP}{{\mathbb{P}}}

 \newcommand{\D}{{\mathcal{D}}}

\renewcommand{\H}{{\mathcal{H}}}
 
 \newcommand{\J}{{\mathcal{J}}}
 \newcommand{\K}{{\mathcal{K}}}

\renewcommand{\O}{{\mathcal{O}}}
\renewcommand{\P}{{\mathcal{P}}}

 \newcommand{\V}{{\mathcal{V}}}
 \newcommand{\W}{{\mathcal{W}}}

\newcommand{\upchi}{{\raise.35ex\hbox{$\chi$}}}
\newcommand{\fD}{{\mathfrak{D}}}

\newcommand{\fH}{{\mathfrak{H}}}

\newcommand{\qand}{\quad\text{and}\quad}

\newcommand{\qfor}{\quad\text{for}\quad}

\newcommand{\qforal}{\quad\text{for all}\quad}

\newcommand{\ran}{\operatorname{Ran}}
\newcommand{\rank}{\operatorname{rank}}

\newcommand{\spn}{\operatorname{span}}

\newcommand{\fngee}{\bbF^+\!(G)}

\newcommand{\sumin}{\sum_{i=1}^n}

\newcommand{\rowt}{(T_1, \ldots, T_n)}
\newcommand{\rows}{(S_1, \ldots, S_n)}

\newcommand{\opL}{\operatorname{{\mathbf L}}}
\newcommand{\mindil}{\operatorname{MinDil}}

\begin{document}

\title[Partially Isometric Dilations of Noncommuting $N$-tuples]%
{Partially Isometric Dilations of Noncommuting $N$-tuples of Operators}
%
\author[M.T.Jury and D.W.Kribs]{Michael~T.~Jury$^1$ and 
David~W.~Kribs$^2$} 
\thanks{{\it key words and phrases.} Hilbert space, operator, row contraction, partial isometry,
minimal dilation, directed graph.} 
\thanks{$^1$first author partially 
supported by a  VIGRE Post-doctoral  Fellowship.} 
\thanks{$^2$second author partially supported by an NSERC
research grant.} 
\address{Department of Mathematics, Purdue University, West Lafayette, 
IN, USA 47907}
\email{jury@math.purdue.edu}
\address{Department of Mathematics and Statistics, University of Guelph, Guelph, ON,
CANADA   N1G 2W1}
\email{kribs@math.purdue.edu}

\date{}
\begin{abstract}
Given a row contraction of operators on Hilbert space and a family of projections on the space
which stabilize the operators, we show there is a unique minimal joint dilation to a  row
contraction of partial isometries which satisfy natural relations. 
For a fixed row contraction the set of all  dilations forms a partially 
ordered set with a largest and smallest element. A
key technical device in our analysis is a connection with directed graphs. We use a Wold
Decomposition for partial isometries to describe the models for these 
dilations,  and discuss how the basic properties of a dilation depend on 
the  row
contraction.  
\end{abstract} 
\maketitle

\section{Introduction}\label{S:intro}

Dilation theory has played a central role in operator theory since 
Sz.-Nagy \cite{SF} proved  in 1953 that
every contraction operator on Hilbert space has a unique minimal  dilation to an
isometry on a larger space. This result was extended to the noncommutative multivariable
setting by Frazho \cite{Fra1} (for $n=2$), Bunce \cite{Bun} (for $2\leq n < \infty$), and
Popescu \cite{Pop_diln} (for $n=\infty$ and uniqueness in general). Specifically, every row
contraction of $n$ operators on Hilbert space was shown to have a joint minimal 
dilation to $n$ isometries on a larger space with mutually orthogonal ranges. The study of
isometries with orthogonal ranges has provided the technical underpinning for a number of far
reaching enquiries (see
\cite{BV1,BJit,DKP,DKS,Jorgen_min,Kribs_curv,Pop_curv} for examples from different
perspectives). While the Frazho-Bunce-Popescu (FBP) dilation  has played a 
role in many of
these instances, there are deep  reasons from the representation theory of infinite dimensional
operator algebras which suggest it may have  limited utility.   

In this paper, we present a dilation theory for
$n$-tuples of operators on Hilbert space which is,
in general, more in tune with properties of the $n$-tuple as compared to the                 
FBP dilation. Given a row contraction $T = \rowt$
and a family of projections which stabilize the operators in a certain sense (such families always
exist), we show there is a unique minimal joint dilation of $T$ to an $n$-tuple of partial
isometries $S = \rows$ which satisfy natural relations. 
This dilation theorem may be regarded as a refinement of a special case of the recent         
Muhly-Solel \cite{MS2} theorem for the more abstract setting of tensor algebras over             
$\ca$-correspondences. 

For fixed $T$, the set of all dilations forms a partially ordered set of directed graphs with a
largest and smallest
element. The smallest element is the Sz.-Nagy dilation (for $n=1$) or the 
FBP dilation (for $n\geq 2$), which corresponds to the directed graph with a
single vertex and $n$ loop edges. There is a `finest'
dilation which is the largest element in the ordering. This dilation is maximal amongst the set of 
all  minimal dilations of $T$ in the sense that if we are given a minimal dilation of $T$, the
corresponding directed graph
is a `deformation' of the graph for the finest dilation.

The main technical drawback of the FBP dilation is that the analogue of the unitary part in the
Wold
Decomposition \cite{Pop_diln} determines a representation of the Cuntz 
algebra $\O_n$, which is an `NGCR'
algebra \cite{Glimm}, and hence its representations cannot be classified up to unitary
equivalence. While this has been accomplished for special classes of row contractions
\cite{DKS}, in general it is not possible. 
On the other hand, for many row contractions, the dilation theory developed here avoids this
problem. Indeed, the analogue of the unitary part here determines a representation of a      
Cuntz-Krieger graph $\ca$-algebra \cite{Kumjian1,Kumjian2}, and there are many graphs for
which the representation theory of the algebra is type I. In fact, Ephrem 
\cite{Eph} has recently
obtained a complete graph-theoretic characterization of when this happens.  

In the first section we discuss the models for this dilation theory and recall the Wold
Decomposition from \cite{JuryK} for families of partial isometries. 
We prove the dilation theorem in the second section, and show how basic
properties of a dilation depend on the dilated row contraction. 
We also describe the class of row contractions for which this dilation theory gives an
improvement on the FBP dilation theory. 
In the final section we discuss 
the partially ordered set of minimal dilations generated by a given row contraction. 

Throughout the paper $n$ is a positive integer or $n=\infty$, but we behave as though $n$ is
finite.

\section{Wold Decomposition}\label{S:wold}

The models for the dilation theory  presented here are $n$-tuples $S = 
\rows$ of (non-zero) 
operators acting on a Hilbert space $\K$ which satisfy the following relations:
\[
(\dagger) \left\{ 
\begin{array}{cl}
(1) & \forall \,1\leq i \leq n, \,\,(S_i^* S_i)^2 = S_i^*S_i \\
(2) & \sum_{i=1}^n S_i S_i^* \leq I \\
(3) & \mbox{$\forall \, 1\leq i,j \leq n,\,\,(S_i^* S_i)(S_j^* S_j) = 0$ or $S_i^* S_i = S_j^*
S_j$} \\
(4) & \mbox{$\forall \, i \, \exists \, j$ such that   $S_i S_i^* \leq S_j^* S_j$} \\
(5) & \mbox{If $\{Q_k\}$ are the distinct elements from $\{S_i^* S_i\}$,} \\ 
   & \mbox{then $\sum_k Q_k = I$.}  
\end{array} 
\right.
\]

Such an $n$-tuple consists of partial isometries with mutually orthogonal ranges, with initial
projections equal or orthogonal, with each final
projection supported by some initial projection, and
distinct initial projections summing to the identity operator. Observe there is a natural directed
graph $G$ (with no sinks) associated with each $n$-tuple which 
satisfies $(\dagger)$. The vertex set $V(G)$ for
$G$ is identified with the index set for $\{Q_k \}_k$, and the  edge set $E(G)$ includes a
directed  edge for each $S_i$; specifically, $S_i$ determines an edge in $G$ from vertex $k$ to
vertex $l$ where $S_i^* S_i = Q_k$ and $S_i S_i^* \leq Q_l$. We will use the orderings of $S =
\rows$ induced by $G$, and write $S = (S_e)_{e\in E(G)}$ when an ordering has been chosen. 

If $S = \rows$ satisfies $(\dagger)$ with $SS^* = \sum_{i=1}^n S_i S_i^* = I$, then we say $S$
is {\it fully coisometric}. From the operator algebra perspective, fully coisometric $n$-tuples
generate what are sometimes called Cuntz-Krieger directed graph 
$\ca$-algebras (see
\cite{Eph,Kumjian1,Kumjian2} for instance).  

At the other extreme, we say $S = (S_e)_{e\in E(G)}$ satisfying $(\dagger)$ is {\it pure} if
\begin{eqnarray}\label{puredilncondition}
\lim_{d\rightarrow \infty}\Big( \sum_{w\in\fngee; \,\, |w|\geq d} 
|| w(S)^* \xi ||^2\Big) = 0  \qforal \xi\in\K.
\end{eqnarray}
Here we denote the {\it semigroupoid} of $G$ by $\fngee$. This is the set of all vertices in $G$
and all finite paths $w$ in the edges $e$ of  $E(G)$, with the natural operations of concatenation
of allowable paths. We write $|w|$ for the number of edges which make up the path $w$, and  
put $w = k_2wk_1$ when the initial and final vertices of $w$ are, respectively, $k_1$ and
$k_2$. The
notation $w(S)$ stands for the partial isometry given by the product $w(S) = S_{e_{i_1}} \cdots
S_{e_{i_m}}$ when $w = e_{i_1} \cdots e_{i_m}$ belongs to $\fngee$.

We now discuss the fundamental examples for the pure case. Let $G$
be a countable directed graph  and let $\K_G = \ell^2(G)$ be the Hilbert space 
with orthonormal basis $\{ \xi_w : w\in \fngee \}$. Define partial isometries on $\K_G$ by 
\[
{\mathbf \opL}_e    \xi_w = \left\{ 
\begin{array}{cl}
\xi_{ew} & \mbox{if $ew\in\fngee$} \\
0 & \mbox{otherwise.}
\end{array}\right.
\]
The operators ${\mathbf \opL}_G = ({\mathbf \opL}_e)_{e\in E(G)}$ are easily seen to be pure
and satisfy $(\dagger)$.
This generalized `Fock space' construction was introduced by Muhly 
\cite{M1} and there is now a growing literature for the 
nonselfadjoint operator algebras generated by such tuples 
\cite{JP,JuryK,KP2,KP1,MS1,MS2}.  The  
$\ca$-algebra generated by a tuple $\opL_G$ is said to be of 
`Cuntz-Krieger-Toeplitz' type since it is the extension of a 
Cuntz-Krieger algebra by the compact operators. 

Every pure tuple $S = (S_e)_{e\in E(G)}$ which satisfies $(\dagger)$ for $G$ is determined by
$\opL_G$ in the following sense:
Let $\V_k$, $k\in V(G)$, be the subspace of $\K_G$ generated by basis vectors from paths
which
begin at vertex $k$, that is, 
$\V_k = \spn\{ \xi_w : w=wk \in \fngee\}$. Then there is a joint unitary
equivalence such that 
\[
S_e \simeq \sum_{k\in V(G)} \! \oplus {\mathbf
\opL}_e^{(\alpha_k)}\Big|_{\V_k^{(\alpha_k)}} 
\qfor e\in E(G), 
\]
where $\alpha_k = \dim \big[Q_k \big( I - \sum_e S_e S_e^* \big) \big]$ 
and recall $\{ Q_k \}_{k\in
V(G)}$ are the distinct projections amongst $\{ S_e^* S_e : e\in E(G) \}$. 
The basic idea is as follows. 
Let $\W = \ran \big( I - \sum_e S_e S_e^* \big)$ be the {\it wandering subspace} \cite{JuryK}
for $S$. A unitary producing the joint equivalence is defined by making a 
natural
identification between  orthonormal bases for the non-zero subspaces of 
the form $w(S) Q_k \W$ and 
corresponding subspaces of $\V_k^{(\alpha_k)}$. 
We refer to the $\alpha_k$ as the {\it vertex multiplicities} in this decomposition. 

The following
Wold Decomposition was established in \cite{JuryK} for $n$-tuples satisfying $(\dagger)$.

\begin{thm}\label{wold}
Let $S=\rows$ be operators on $\K$ satisfying $(\dagger)$ and let $S = (S_e)_{e\in E(G)}$ be
an induced ordering. Then these operators 
are jointly unitarily equivalent to the direct sum of a pure $n$-tuple and a
fully coisometric $n$-tuple which both satisfy $(\dagger)$ for the directed graph $G$. In other
words, there is a unitary $U$ and a fully coisometric $n$-tuple  $(V_e)_{e\in E(G)}$ such that 
\begin{eqnarray}\label{woldid}
\,\,\,\,\,\,\,\,\, US_eU^* = 
 V_e \oplus \Big( \sum_{k\in V(G)}\oplus \opL_e^{(\alpha_k)}\Big|_{\V_k^{(\alpha_k)}}
\Big) \qfor e\in E(G), 
\end{eqnarray}
and the $\alpha_k$ are determined as above. 

Let $\K_p = \sum_{w\in\fngee} \oplus w(S) \W$ where $\W = \ran \big( I - \sum_e S_e S_e^*
\big)$ and let $\K_c = ( \K_p)^\perp$. 
The subspaces $\K_c$ and $\K_p$ reduce $S = (S_e)_{e\in E(G)}$, and the restrictions
$S_e|_{\K_c}$ and 
$S_e|_{\K_p}$ determine the joint unitary equivalence in $(\ref{woldid})$. 
This decomposition is unique in the sense that if $\V$ is a
subspace of $\K$ which reduces $S = (S_e)_{e\in E(G)}$, and if the restrictions $\{ S_e|_\V :
e\in E(G)\}$
are pure, respectively fully coisometric, then $\V \subseteq \K_p$, respectively $\V \subseteq
\K_c$. 
\end{thm}

We finish this section by identifying a large class of pure row 
contractions which will
be used in the sequel. 
Let $\fH = \{\H_k :k\in \J\}$ be a countable collection of Hilbert spaces.  
Let $G$ be a countable
directed graph with vertex
set $V(G) = \J$. We define $\ell^2(G,\fH)$ to be the Hilbert space given by the
$\ell^2$-direct
sum $\ell^2(G,\fH) = \sum_{w\in \bbF^+(G)} \oplus \H_w$ where $\H_w \equiv  \H_k$ when
$w
= wk$, that is, the initial vertex of $w$ is $k$. For each non-zero $\H_k$ choose an orthonormal
basis $\{ \xi_{j}^{(k)} 
\}$, and for $w=wk \in \bbF^+(G)$ let $\{ \xi_{j}^{(w)} \}$ be the corresponding orthonormal
basis for the $w$th coordinate space $\H_w$ of $\ell^2(G,\fH)$. 
Then the {\it canonical (pure) shift} on $\ell^2(G,\fH)$ consists of operators  $(L_e)_{e\in
E(G)}$ defined on $\ell^2(G,\fH)$ by  
\[
L_e \xi_{j}^{(w)} = 
\left\{ \begin{array}{cl}
\xi_{j}^{(ew)} & \mbox{if $ew\in\fngee$} \\
0 & \mbox{otherwise.}
\end{array}\right.
\]
It is easy to see that every canonical shift $(L_e)_{e\in E(G)}$ is pure and satisfies $(\dagger)$,
and hence
Theorem~\ref{wold} explicitly gives its form up to joint unitary 
equivalence. In particular, the fully coisometric part is vacuous and 
the vertex multiplicities are given by $\alpha_k = \dim(\H_k)$ for $k\in 
V(G)$.

\section{Minimal Partially Isometric Dilations}\label{S:dilation}

Let $T = \rowt$ be operators on a Hilbert space $\H$ such that $TT^* = \sum_{i=1}^n T_i
T_i^* \leq I_\H$. 
We say an $n$-tuple $S = \rows$ of operators on a Hilbert space $\K \supseteq
\H$ is a {\it minimal partially isometric dilation} of $T$ if the following conditions hold:
\begin{itemize}
\item[$(i)$]
$S = \rows$ satisfy the relations $(\dagger)$. 
\item[$(ii)$]
$\H$ reduces each $S_i^* S_i$, $1\leq i \leq n$, and $\H$ is invariant for each $S_i^*$ with
$S_i^*|_{\H} = T_i^*$, $1\leq i \leq n$. 
\item[$(iii)$]
 $\K = \H \vee \Big( \bigvee_{i_1,\ldots,i_k;\, k\geq 1} S_{i_1}\cdots S_{i_k} \H \Big)$. 
\end{itemize}

Given $T = \rowt$, consider all countable families  $\P = \{P_k :k\in\J\}$ 
of projections on $\H$
which {\it stabilize} $T$ in the following sense:
\begin{eqnarray}\label{projnident}
P_k T_i, \,\, T_i P_k \,\in \{T_i, 0\}, \,\, 1\leq i \leq n,
\qand \sum_{k\in\J} P_k = I_\H. 
\end{eqnarray}
For each $i$, it will be convenient to let $k_s = s(T_i)$ and $k_r = r(T_i)$ be the elements of
$\J$ such that $T_i P_{k_s} = T_i$ and $P_{k_r} T_i = T_i$. 
To avoid pathologies we shall assume each $T_i$ is non-zero. (If some $T_i 
=0$, there is ambiguity in the choice of $k_s, k_r$.)
Further, we clearly lose no generality in restricting our attention to 
families $\P$ such that there is no $P_k$ with $T_i P_k = 0$ for all $i$. 
We show there is a minimal dilation of $T$ generated by each such family of projections.

Given a family $\P = \{P_k \}_{k\in\J}$ which satisfies 
(\ref{projnident}), we let  $I_\P$ be the
projection on the Hilbert space direct sum $\H^{(n)}$ defined by the $n\times n$
diagonal matrix with $(i,i)$ entry equal to $P_{k_i}$ where $k_i = s(T_i)$. Observe that
the relations (\ref{projnident}) guarantee that $I_\P - T^* T\geq 0$ is a 
positive operator on
$\H^{(n)}$, here regarding $T$ as a row matrix. 
Thus we may define a {\it defect operator} for $T = \rowt$ on $\H^{(n)}$ by $D\equiv
D_{\P,T} = \big( I_\P - T^* T \big)^{1 / 2}$. Let $\D = \ol{D \H^{(n)}}$.
Further let $E_i : \H \rightarrow \H^{(n)}$ be the injection of $\H$
onto the $i$th coordinate space of $\H^{(n)}$ for $1\leq i \leq n$. Consider
the operators $D_i = DE_i : \H \rightarrow \D $ for $1\leq i \leq n$. 

\begin{lem}\label{mainlemma}
If $r(T_i) \neq r(T_j)$, then the range subspaces $\ran (D_i)$ and $\ran (D_j)$ are orthogonal. 
\end{lem}

\Prf
It suffices to show that $\ran (D^{2a} E_i)$ and $\ran (D^{2b} E_j)$ are 
orthogonal for $a,b\geq 1$; then a
standard functional calculus argument can be applied. 
Recall that $D^2 = I_\P - T^*T$ is an $n\times n$ matrix which acts on $\H^{(n)}$. The
operator $D^2 E_i$ picks out the $i$th column of $D^2$. When $r(T_i) \neq r(T_j)$, it follows
from the identities (\ref{projnident}) that in the $i$th and $j$th columns of $D^2$ there are no
rows $m$ such that both the $(m,i)$ and $(m,j)$ entries are non-zero. 
This property is easily seen to carry over to the self-adjoint powers 
$(D^2)^a$. 
Hence  
$D^{2a}E_i$ and $D^{2b}E_j$ have  orthogonal ranges for $a,b \geq 1$ when 
$r(T_i) \neq r(T_j)$. 
\bx

We will use these operators to define generalized Schaffer
matrices \cite{Fra2,MS2,Pop_diln,Schaffer} in the following proof. 

\begin{thm}\label{mindiln}
Let $T = \rowt$ be operators on a Hilbert space $\H$ such that $\sum_{i=1}^n T_i T_i^* \leq
I_\H$. Let $\P = \{P_k \}_{k\in\J}$ be a family of projections which stabilize $T$ as in 
$(\ref{projnident})$.   Then there is a minimal partially isometric dilation $S = \rows$ of $T$ on
a Hilbert space $\K \supseteq \H$ with $\P = \{ S_i^* S_i|_{\H}: 1\leq i \leq n \}$. This dilation
is unique up to joint unitary equivalence which fixes $\H$. 
\end{thm}

\Prf
By Lemma~\ref{mainlemma} we may decompose $\D$ into the orthogonal direct sum
$\D = \sum_{k\in V(G)} \oplus \D_k$, where $\D_k = \bigvee_{r(T_i) = k} 
\ol{D_i \H}$ are
subspaces of $\D$ and $G$ is the directed graph (with no sinks) 
determined by $\P$ and the relations (\ref{projnident}). Put $\fD = \{ \D_k \}_{k\in V(G)}$.  
Let $\K$ be the Hilbert space $\K = \H \oplus \ell^2(G,\fD)$ and let $\H$ and $\ell^2(G,\fD)$ be
embedded
into $\K$ in the natural way. For the rest of this proof it is convenient to re-label $T = \rowt$ as
$T =
(T_e)_{e\in E(G)}$ by using a natural ordering induced by (\ref{projnident}). We
shall carry this notation over to the operators $D_i, E_i$, denoting them by $D_e, E_e$. 
For each $e\in E(G)$ define operators $S_e : \K \rightarrow \K$ by 
\[
S_e = \big( T_e + D_e \big) \oplus L_e, 
\]
where $(L_e)_{e\in E(G)}$ is the canonical shift on $\ell^2(G, \fD)$. 

We first verify $(i)$ and $(ii)$ for a minimal dilation. 
Observe that 
\[
T_e^* T_f + D_e^* D_f = T_e^* T_f + E_e^* \big( I_\P - T^* T \big) E_f. 
\]
When $e=f$, this identity yields $T_e^* T_e + D_e^* D_e = P_{k_0}$ where $k_0 = s(T_e)$.
On the other hand, if $e\neq f$, then 
$T_e^* T_f + D_e^* D_f = T_e^* T_f - T_e^* T_f = 0$. 
It follows that the operators $S_e$ are partial isometries with $T_e = P_\H S_e|_{\H} = \big(
S_e^*|_{\H} \big)^*$, and initial projections which satisfy $\{ S_e^* S_e|_{\H}: e\in E(G) \} =
\{ P_k :k\in V(G)\}$.
Moreover, the ranges of the $S_e$ are mutually orthogonal, $S_e^* S_f = 0$ for $e\neq f$, and
hence $\sum_e S_e S_e^* \leq I_\K$. 
Lastly, by construction each range projection $S_e S_e^*$ is supported by an initial
projection, $S_e S_e^* \leq S_f^* S_f$ for some $f\in E(G)$, and the 
distinct initial projections $\{ Q_k \}_{k\in V(G)} = \{ S_e^* S_e 
\}_{e\in E(G)}$ sum to the identity. 

To verify minimality, first notice that $\H \bigvee S_e \H = \H \oplus \D$. But 
\begin{eqnarray*}                                      
\K \ominus (\H \oplus \D) &=& \ell^2(G,\fD) \ominus \D \\ &=& \sum_{e\in E(G)} \!\!\!\!
\oplus\,\, S_e ( \ell^2(G,\fD))  =
\sum_{w\in\fngee; \, |w|\geq 1}\!\!\!\!\!\!\!\!\!\! \oplus\,\,\, w(S) \D, 
\end{eqnarray*}
and thus we have $\K = \H \vee \Big( \bigvee_{w\in\fngee;\,|w|\geq 1} w(S) 
\H\Big)$.

Finally, the uniqueness assertion is that if $S^\prime = (S_1^\prime, \ldots, S_n^\prime)$ on
$\K^\prime \supseteq \H$ is another minimal dilation of $T$ with respect to $\P$, then there is a
unitary $U : \K \rightarrow \K^\prime$ such that $U|_\H = I_\H$ and $U^* S_i^\prime U = S_i$
for $1\leq i \leq n$. This proof is a relatively simple adaptation of the single variable case
\cite{SF}, hence we omit the details. 
\bx

\begin{rem}
In the case that the family $\P = \{I \}$ is a singleton, Theorem~\ref{mindiln} collapses to the
Sz.-Nagy dilation theorem \cite{SF} when $n=1$ and the FBP dilation 
theorem 
\cite{Fra1,Bun,Pop_diln} when $2\leq n \leq \infty$. This is the only case for which the minimal
dilation consists entirely of isometries.    
In its most general form, Theorem~\ref{mindiln} may be regarded as a refinement of the         
Muhly-Solel dilation theorem for a subclass of the representations considered in \cite{MS2}. In
the language of \cite{MS2}, a row contraction $T$ and
collection of projections $\P$ satisfying (\ref{projnident}) can be seen to induce a {\it covariant
representation} of a $\ca$-correspondence generated by $T$ and $\P$. These representations
form a subclass of those considered in \cite{MS2}, and the class of all such representations are
shown to have minimal dilations. Hence the basic existence of minimal dilations in our
setting can be
deduced from \cite{MS2}. However, our short spatial proof and the particular details we obtain
are not easily seen there.  Furthermore, we suggest that the results of the current paper provide a
more accessible dilation theory for row contractions, as the abstract
machinery of Hilbert
modules and $\ca$-correspondences is not required in the formulation here. 
\end{rem}

We next discuss how properties of $T$ can be used to identify properties of its minimal
dilations. 

\begin{prop}\label{pureident}
Every minimal partially isometric dilation of $T = (T_e)_{e\in E(G)}$ is pure if and only if 
\begin{eqnarray}\label{purecondition}
\lim_{d\rightarrow \infty}\Big( \sum_{w\in\fngee; \, |w|\geq d} 
|| w(T)^* \xi ||^2 \Big) = 0  \qforal \xi\in\H.
\end{eqnarray}
\end{prop}

\Prf
If $S = (S_e)_{e\in E(G)}$ is a  pure minimal dilation of $T$, then 
$S_e^*|_\H = T_e^*$ for
$e\in E(G)$ and (\ref{purecondition}) follows from
the corresponding identity (\ref{puredilncondition}) for $S$. 
Conversely, when (\ref{purecondition}) holds we may use the $(\dagger)$ relations to obtain the
necessary estimates which  show that
(\ref{puredilncondition}) holds for every minimal dilation $S$ of $T$. 
\bx

\begin{cor}
If $\sumin T_i T_i^* \leq r I$, with $r < 1$, then every minimal partially isometric dilation of $T
=\rowt$ is pure. 
\end{cor}

Next we obtain detailed information on the pure part of a dilation.

\begin{prop}\label{vertexmult}
Let $S = (S_e)_{e\in E(G)}$ be a
minimal partially isometric dilation of $T = (T_e )_{e\in E(G)}$ with respect to the projections
$\P = \{P_k \}_{k\in V(G)}$. Then for $k\in V(G)$ we have 
\begin{eqnarray}\label{vertexrank}
\rank \Big( Q_k \big( I_\K - \sum_e S_e S_e^* \big) \Big) = 
\rank \Big( P_k \big( I_\H - \sum_e T_e T_e^* \big) \Big).
\end{eqnarray} 
\end{prop}

\Prf
Recall $\{Q_k \} = \{S_e^* S_e \}$. By the Wold Decomposition and
Theorem~\ref{mindiln} we may assume that $Q_k|_\H = P_k$ for $k\in V(G)$. 
Fix $k\in V(G)$.
Observe the $(\dagger)$ relations imply $Q_k$ commutes with $P = I_\K - \sum_e S_e S_e^*$.
Let $R_k$ be the projection $R_k = Q_k P$ and let $P_\H$ be the projection of $\K$ onto $\H$.
The minimality of the dilation ensures the subspace $P\K$ does not intersect $\H^\perp$, and
hence neither does the subspace $Q_k P \K = P Q_k \K$. Thus $P_\H (Q_k P) 
P_\H$ has the same
rank as $Q_k P$ (even though $Q_k P \K$ is not contained in $\H$ in general). But notice that 
\[
P_\H Q_k P \Big|_\H = P_\H \Big( Q_k \big( I_\K - \sum_e S_e S_e^* \big) \Big)P_\H
\Big|_\H 
= P_k \big( I_\H - \sum_e T_e T_e^* \big),  
\]
and the result follows. 
\bx

\begin{cor}\label{coisometric}
Every minimal partially isometric dilation of $T = \rowt$ is fully coisometric if and only if 
$
\sumin T_i T_i^* = I_\H.
$
\end{cor}

\begin{rem}
The identity (\ref{vertexrank}) shows how to compute the vertex multiplicities for a minimal
dilation strictly in terms of the dilated row contraction $T$ and the projection family $\P$. Thus, 
by Theorem~\ref{wold} this gives a method for explicitly finding the pure part of a dilation. 

On the other hand, Corollary~\ref{coisometric} identifies when the fully 
coisometric case occurs in terms of $T$. 
In complete generality it is not possible to explicitly describe the fully 
coisometric part of a minimal
dilation. 
As mentioned above, the representation theory of $\O_n$ is the obstacle. However, note that the
fully coisometric part of a minimal dilation here will determine a representation of a 
Cuntz-Krieger directed graph $\ca$-algebra $\ca (G)$. Ephrem \cite{Eph} characterizes when
the representation theory of such algebras is type~I strictly in terms of the directed graph $G$.
Interestingly, in the case of finite graphs his graph-theoretic condition can be seen to be precisely
the condition obtained by the
second author and Power \cite{KP2,KP1} as a description of when a nonselfadjoint
`free semigroupoid  algebra' is partly free. Specifically, $\ca (G)$ is 
type~I if and only if 
the following two conditions hold:
\begin{itemize}
\item[$(i)$]
$G$ contains no double-cycles; there are no distinct cycles   $w_1 
=xw_1x$,
$w_2 =xw_2x$  at a vertex $x$ in $G$.
\item[$(ii)$]
Given a non-overlapping infinite directed path in $G$, there are only finitely many ways to exit
and return to the path. 
\end{itemize}

Thus, whenever the $G$ obtained in a minimal dilation of $T$ satisfies these conditions, the
dilation theory here gives an improvement on the dilation theory derived from the FBP dilation.
For example, let $\H$ be a Hilbert space and let $T_1,T_2,T_3$ be operators on $\H$ such that
$T_1$ is a co-isometry and $(T_2, T_3)$ forms a row contraction with $T_2T_2^* + T_3T_3^*
= I$.
Define a row contraction $V = (V_1,V_2,V_3)$ on $\H \oplus \H$ by 
\[
V_1 = 
\left[
\begin{matrix}
T_1 & 0 \\
0& 0 
\end{matrix}
\right], \quad\quad 
V_2 = \left[
\begin{matrix}
0 & 0 \\
0& T_2 
\end{matrix}
\right], \quad\quad 
V_3 = 
\left[
\begin{matrix}
0 & 0 \\
T_3& 0 
\end{matrix}\right]. 
\]
Let $P_i\equiv P_\H$, $i=1,2$, be the projections of the  direct sum $\H 
\oplus \H$ onto its two
coordinate spaces.
Observe that $VV^* = \sum_{i=1}^3 V_i V_i^* = I$. Hence, the minimal partially isometric
dilation of $V$ with respect to $\P = \{P_1, P_2\}$ determines a representation of the 
$\ca$-algebra $\ca(G)$ where $G$ is the directed graph with two vertices, a loop edge over each
vertex, and a directed edge from the first to the second vertex. As $G$ satisfies the above
conditions, $\ca(G)$ is type~I and hence GCR \cite{Glimm}.  
That is, every representation of $\ca(G)$ can be obtained as a direct 
integral of irreducible subrepresentations \cite{Arvinvite}. 
\end{rem}

%

\section{Partial Orders on Minimal Dilations}\label{S:poset}

%
%

Let $T = \rowt$ be a fixed row contraction on $\H$. As a convenience, in 
this section we assume that 
$\H = \vee_i (\ran(T_i) \vee \ran(T_i^*))$. (Observe that the restrictions 
of each $T_i$ and $T_i^*$ to the orthogonal complement of this joint 
reducing subspace are zero.) 

\begin{lem}\label{commute}
If $\P_1$ and $\P_2$ are families of projections which stabilize $T$ as in 
$(3)$, then the families $\P_1$ and $\P_2$ are mutually commuting. 
\end{lem}

\Prf
Let $P_\alpha \in \P_1$ and $P_\beta \in \P_2$. Then for $1\leq i \leq n$, 
\[
P_\alpha P_\beta T_i = P_\beta P_\alpha T_i = 
\left\{ \begin{array}{cl}
T_i & \mbox{if $P_\beta T_i = T_i = P_\alpha T_i$} \\
0 & \mbox{if $P_\beta T_i =0$ or $P_\alpha T_i = 0$.}
\end{array}\right.
\]
Similarly, $T_i P_\alpha P_\beta = T_i P_\beta P_\alpha$, so that 
$P_\beta P_\alpha T_i^* = P_\alpha P_\beta T_i^*$. 
By the assumption above, the ranges of $\{T_i,T_i^* : 1\leq i \leq n\}$ 
are dense inside $\H$, hence the result follows. 
\bx

Let $G_\P$ be the directed graph determined by the relations (\ref{projnident}) for a family of
projections $\P$ which stabilize $T$. 
Using the Wold Decomposition and Theorem~\ref{mindiln}, we may identify 
the set of all graphs $G_\P$ with the set
$\mindil(T)$ of all equivalence classes of minimal partially isometric dilations of $T$.  
Define a partial ordering on $\mindil(T)$ by: 
 $G_{\P_1} \leq G_{\P_2}$ if and only if  
\[
\forall\, P\in\P_1 \quad \exists \, \bbP \subseteq \P_2 \quad\mbox{such that}\quad    
P = \sum_{P_\alpha \in \bbP} P_\alpha.
\]
This set has a natural join operation defined by $G_{\P_1} 
\bigvee G_{\P_2} \equiv
G_{\P_1 \bigvee \P_2}$ where  
\[
\P_1 \bigvee \P_2 = \Big\{ P_1 \wedge P_2 : P_i \in \P_i, \, i=1,2 \Big\}, 
\]
and $P_1 \wedge P_2 = P_1 P_2 = P_2 P_1$ is the projection onto the 
intersection of the range subspaces for $P_1,
P_2$ by Lemma~\ref{commute}.

In terms of the directed graph structures, the relation $G_{\P_1} \leq G_{\P_2}$ means
$G_{\P_2}$ may be {\it deformed},
by identifying certain vertices in $V(G_{\P_2})$, to obtain $G_{\P_1}$.
Conversely, to every deformation of $G_{\P_2}$ there corresponds an element of $\mindil(T)$. 

\begin{prop}\label{smallest}
The partially ordered set $\mindil(T)$ has a largest element and a smallest element. 
\end{prop}

\Prf
The smallest element of $\mindil(T)$ is clearly the minimal {\it isometric} dilation of $T$
\cite{SF,Fra1,Bun,Pop_diln},
corresponding to the directed graph with a single vertex and  $n$ distinct loop edges. On the
other hand, if we let $\ol{\P} = \{ \P_\alpha \}_{\alpha\in\bbA}$ be the set of all sets of
projections which satisfy (\ref{projnident}) for $T$, we may define a largest element of
$\ol{\P}$ by 
$
\P_0 = \big\{ \bigwedge_{\alpha\in\bbA} P_\alpha : P_\alpha \in \P_\alpha \big\}.
$
Indeed, the family $\P_0$ is clearly a set of pairwise orthogonal 
projections which stabilize $T$. Further, by Lemma~\ref{commute} the 
$P_\alpha$ from distinct $\P_\alpha$ commute, and hence $\P_0$ determines 
a partition of the identity which is the supremum of  $\ol{\P}$.
\bx

The unique minimal partially isometric dilation $S_0$ which corresponds to
$\P_0$  may be regarded as the `finest' of all
the minimal partially isometric dilations of $T = \rowt$. It is the minimal
dilation which  best reflects the joint  behaviour of the $T_i$. 
Of course, in many instances this will be the minimal isometric dilation of $T$, when $\P =
\{I\}$ is the only projection family which stabilizes $T$, but in general 
there may be non-trivial
families of such projections.  
Amongst the set of all 
directed graphs $G_\P$ which come from minimal dilations of $T$, the directed graph
$G_{\P_0}$ will be the largest in the sense that   any other directed
graph in this set can be
obtained from $G_{\P_0}$ by a series of deformations.


\vspace{0.1in}

{\noindent}{\it Acknowledgements.}
This paper was partly motivated by Paul Muhly's talk at the 2003 Great Plains Operator Theory
Symposium. 
We are grateful to Alex Kumjian and Menassie Ephrem for helpful discussions. 
We would also like to thank members of the Department of Mathematics at 
Purdue University for
kind hospitality during the preparation of this article.


%

\end{document}